# Comparative work for the source identification in parabolic inverse problem based on Taylor and Chebyshev wavelet methods

Gopal Priyadarshi[a,b] and Sila Ovgu Korkut[c]
[a] Applied Mathematics and Computational Sciences
King Abdullah University of Science and Technology (KAUST), Saudi Arabia
[b] Department of Mathematics
S.M.D. College, Patliputra University, India
gopal.priyadarshi@kaust.edu.sa

[c] Department of Engineering Sciences
Izmir Katip Celebi University
Izmir, 35640, Turkey
silaovgu@gmail.com

**Abstract**: In this article, we study wavelet collocation methods based on Taylor and Chebyshev wavelets for the source identification in parabolic inverse problem. In the proposed method, highest order derivative is written in terms of Taylor and Chebyshev wavelet series and required unknown terms are obtained using successive integration. Taylor series approximation has been utilized to obtain the source control parameter. Convergence analysis is carried out in order to guarantee the accuracy of the method. Numerical results have been obtained based on the proposed methods and it is shown that Taylor wavelet method provide us better result than the Chebyshev wavelet method. CPU time has also been shown to ensure the efficiency of the method.

## 1 Introduction

In recent years, there has been increasing interest in the development of efficient and accurate numerical methods for parabolic inverse problem with source control parameter. These problems have many applications in thermoelasticity, heat conduction process, control theory, biochemistry and population dynamics.

In this study, we propose collocation methods based on Taylor and Chebyshev wavelets to identify the source control parameter appear in parabolic inverse problem given by

$$\frac{\partial y}{\partial t}(x,t) = A\frac{\partial^2 y}{\partial x^2}(x,t) + B\frac{\partial y}{\partial x}(x,t) + X(t)y(x,t) + \psi(x,t), \ x \in (0,1), \ 0 \leq t \leq T \quad (1.1)$$

with initial condition

$$y(x,0) = y_0(x), \ x \in (0,1) \quad (1.2)$$

and Dirichlet boundary conditions

$$y(0,t) = f_0(t), \qquad y(1,t) = f_1(t), \quad 0 \leq t \leq T, \quad (1.3)$$

subject to the overspecified condition

$$y(x_{in}, t) = Q(t), \quad 0 \leq t \leq T, \tag{1.4}$$

where $A, B$ are constants, $x_{in}$ is a fixed point such that $0 < x_{in} < 1$ and $\psi, y_0, f_0, f_1, Q$ are known and sufficiently smooth functions. It is also assumed that $Q(t) \neq 0$. We have to determine unknown functions $y(x, t)$ and $X(t)$ simultaneously, where $y(x, t)$ stands for the solution of the given problem and $X(t)$ is the source control parameter.

The existence, uniqueness and regularity results for the parabolic inverse problem is discussed by Cannon at al.[2,3] and Prelipko et al.[20,21] Unfortunately, these type of problems turn out to be ill posed in the sense of Hadamard because the problem is not stable. In other words, solution of the problem does not depend continuously on the given data. Various methods e.g. the projective method,[18] dual reciprocity boundary element methods,[9] potential logarithmic method[19] etc. have been developed for the source identification in the parabolic problem.

Numerical methods based on finite difference has been studied by Dehghan et al.[6] They developed various other numerical methods which can be seen in[5,7] and references therein. Liao et al.[17] proposed high-order compact finite difference method to identify the parameter for the parabolic inverse problem. The proposed numerical scheme is fourth order accurate in both time and space. Ritz least square method has been developed by Khorshidi et al.[15] to identify the control parameter in parabolic inverse problem.

In last few years, there is an impressive increase in the use of mesh-free methods in many problems that appear in several areas of science and engineering. Ignoring the mesh production and/or mesh refinement, the meshless methods lead to a significant reduction in the computational cost. One of the most popular technique of mesh-free methods is wavelet method. Wavelet methods have been applied in various areas of science and engineering.[12] Wavelet method based on Legendre multiscaling function has been proposed by Yousefi et al.[28] to solve parabolic inverse problem. Recently Rathish et al.[22] proposed wavelet collocation method based on Haar wavelets to solve multidimensional parabolic inverse problems.

Among various wavelet methods, Chebyshev wavelet method and Taylor wavelet method have great deal of influence in recent years. Wavelet collocation method based on Chebyshev wavelet for the Burger-Huxley equation has been developed by Çelik et al.[4] Farooq et al.[10] developed Chebyshev wavelet method for fractional delay differential equations. Baghani et al.[1] developed second Chebyshev wavelets method for solving finite-time fractional linear quadratic optimal control problems. Köse et al.[8] used Chebyshev wavelet method to solve nonlinear time-fractional Schrödinger equation. More recently, the fourth kind of Chebyshev wavelet method has been utilized in[26] for obtaining a numerical solution of multi-term variable order fractional differential equations. Taylor wavelet method has been proposed in[14] to solve Bratu-type equations. The method has been applied to Lane-Emden equation in both linear and nonlinear sense in[11] and to fractional pantograph equation by Vichitkunakorn et al.[27] A numerical method based on the Taylor wavelet has been proposed in[25] to solve fractional delay differential equations. Rostami et al.[23] developed a numerical method based on the operational matrices of the Taylor wavelet along with the Newton method for solving a class of nonlinear partial integro-differential equations with weakly singular kernels. A numerical method based on Taylor wavelet has been proposed by Korkut et al.[16] to solve general type of KdV-Burgers' equation.

In addition to the above-mentioned valuable studies, the main contribution of this study is to present a comparative work based on Chebyshev and Taylor wavelet methods for the source identification in parabolic inverse problem. A rigorous functional error estimations of these wavelet methods are studied using the important characteristics such as orthogonality of the Chebyshev wavelet and normality of the Taylor wavelet, and etc. In what follows, the convergence analysis is carried out in order to ensure the accuracy of the proposed method.

In this paper, the topics are treated under four headings: any theoretical background of the method including the definitions and implementation of the method has been presented in Section 2. Convergence analysis is carried out in Section 3. Numerical results and discussions are presented in Section 4. The study has been finalized with a brief conclusion in Section 5.

# 2 Method Statements

The aim of this section is to present the whole methodology. To do so, we first describe the Chebyshev and Taylor wavelets and the required first and second integral form of both the wavelets. Secondly, for any interested reader, the implementation of the method bas been described explicitly.

## 2.1 Basic Definitions

Wavelets emerge as a strong and attractive tool which provide the use of quantities expressed by various measures of length in a natural way. Wavelets are comprised by dilation and translation of the mother wavelet by some factors. That's why basically wavelets can vary with respect to the special choices of the mother wavelet. Due to the objective of the present study two types of mother wavelets have been taken into account throughout the study: the Taylor wavelet and the Chebyshev wavelet. We point the reader to our fundamental references for these wavelets,[4,13,14,27]

Before describing the methods, it is vital to note that to alleviate the complexity of the exposition the notations have been classified for their intended use. That is, $I_{nm}(x)$ represents the wavelet methods for either the Taylor wavelets or Chebyshev wavelets. In addition, $R_{nm}(x)$ and $S_{nm}(x)$ are used to denote the first and second integral forms, respectively, for both methods.

Let $n = 1, 2, \ldots, 2^{k-1}$ for $k \in \mathbb{Z}^+$ and $m = 0, 1, \ldots M-1$ where $M$ stands for the degree of the polynomial, the Taylor wavelets defined by Keshavarz et al.[14] can be described as follows:

$$I_{nm}(x) = \begin{cases} 2^{\frac{k-1}{2}} \hat{T}_m \left(2^{k-1}x - n + 1\right), & \text{if } \frac{n-1}{2^{k-1}} \leq x < \frac{n}{2^{k-1}}; \\ 0, & \text{otherwise,} \end{cases} \qquad (2.1)$$

where $\hat{T}_m(x) = \sqrt{2m+1}\, x^m$ for $m = 0, 1, 2, \ldots, M-1$. Here, $x^m$ represents the Taylor polynomial of degree $m$ which forms an orthonormal basis over the interval $[0, 1]$.

Likewise, the Chebyshev wavelets can be identified as follows:

$$I_{nm}(x) = \begin{cases} \gamma_m 2^{\frac{k-1}{2}} C_m \left(2^k x - 2n + 1\right), & \text{if } \frac{n-1}{2^{k-1}} \leq x < \frac{n}{2^{k-1}}; \\ 0, & \text{otherwise,} \end{cases} \quad (2.2)$$

where

$$\gamma_m = \begin{cases} \frac{\sqrt{2}}{\sqrt{\pi}}, & \text{if } m = 0; \\ \frac{2}{\sqrt{\pi}}, & \text{if } m \neq 0, \end{cases} \quad (2.3)$$

for $n = 1, 2, \ldots, 2^{k-1}$, $k \in \mathbb{Z}^+$ and $m = 0, 1, \ldots M - 1$. Notice that $C_m(x)$ stands for the Chebyshev polynomials of the first kind of degree $m$. It is worth emphasizing that Chebyshev wavelets are orthogonal weighted by $\omega\left(2^k x - 2n + 1\right) = \dfrac{1}{\sqrt{1 - \left(2^k x - 2n + 1\right)^2}}$.

In the light of the above-mentioned definitions, any function $f(x) \in L^2_\omega[0,1]$ can be expanded in terms of both the Chebyshev wavelets and the Taylor wavelets as follows:

$$f(x) = \sum_{n=1}^{\infty} \sum_{m=0}^{\infty} d_{nm} I_{nm}(x), \quad (2.4)$$

where $d_{nm}$ are the wavelet coefficients that can vary according to the choice of wavelet methods and they are constituted by

$$d_{nm} = \langle f(x), I_{nm}(x) \rangle = \begin{cases} \displaystyle\int_0^1 f(x) I_{nm}(x) dx, & \text{for the Taylor wavelets}; \\ \displaystyle\int_0^1 f(x) I_{nm}(x) \omega(x) dx, & \text{for the Chebyshev wavelets,} \end{cases} \quad (2.5)$$

Moreover, the approximation of $f(x)$ is given by

$$f(x) \simeq \sum_{n=1}^{2^{k-1}} \sum_{m=0}^{M-1} d_{nm} I_{nm}(x) = \boldsymbol{D}\boldsymbol{I}(x), \quad (2.6)$$

where $\boldsymbol{D}$ is a $1 \times 2^{k-1} M$ vector given by

$$\boldsymbol{D} = \left[d_{10}, d_{11}, \ldots, d_{1(M-1)}, d_{20}, d_{21}, \ldots, d_{2(M-1)}, \ldots, d_{2^{k-1}0}, d_{2^{k-1}1}, \ldots, d_{2^{k-1}(M-1)}\right], \quad (2.7)$$

and

$$\boldsymbol{I}(x) = \left[I_{10}(x) \ldots I_{1(M-1)}(x), I_{20}(x) \ldots I_{2(M-1)}(x), \ldots, I_{2^{k-1}0}(x) \ldots I_{2^{k-1}(M-1)}(x)\right]^T. \quad (2.8)$$

Note that, by replacing the suitable collocation points into Eq. (2.8) $\boldsymbol{I}(x)$ returns $2^{k-1} M \times 2^{k-1} M$ matrix.

We conclude this introductory section with the required integral forms of both wavelets. The first and the second integral forms of the Taylor wavelet is:

$$R_{nm}(x) = \begin{cases} 0, & \text{if } 0 \leq x < \frac{n-1}{2^{k-1}}; \\ \dfrac{2^{\left(m+\frac{1}{2}\right)(k-1)} m! \sqrt{2m+1}}{(m+1)!} \left(x - \frac{n-1}{2^{k-1}}\right)^{m+1}, & \text{if } \frac{n-1}{2^{k-1}} \leq x < \frac{n}{2^{k-1}}; \\ \dfrac{2^{\left(m+\frac{1}{2}\right)(k-1)} m! \sqrt{2m+1}}{(m+1)!} \left(x - \frac{n-1}{2^{k-1}}\right)^{m+1} - P_1(x), & \text{if } \frac{n}{2^{k-1}} \leq x \leq 1, \end{cases} \quad (2.9)$$

and

$$S_{nm}(x) = \begin{cases} 0, & \text{if } 0 \leq x < \frac{n-1}{2^{k-1}}; \\ \frac{2^{(m+\frac{1}{2})(k-1)} m! \sqrt{2m+1}}{(m+2)!} \left(x - \frac{n-1}{2^{k-1}}\right)^{m+2}, & \text{if } \frac{n-1}{2^{k-1}} \leq x < \frac{n}{2^{k-1}}; \\ \frac{2^{(m+\frac{1}{2})(k-1)} m! \sqrt{2m+1}}{(m+2)!} \left(x - \frac{n-1}{2^{k-1}}\right)^{m+2} - P_2(x), & \text{if } \frac{n}{2^{k-1}} \leq x \leq 1, \end{cases} \quad (2.10)$$

where

$$P_i(x) = \sum_{j=0}^{m} \binom{m}{j} \frac{2^{(j+\frac{1}{2})(k-1)} j! \sqrt{2m+1}}{(j+i)!} \left(x - \frac{n}{2^{k-1}}\right)^{j+i}, \quad i = 1, 2. \quad (2.11)$$

Furthermore, the first and the second integral forms of the Chebyshev wavelets[4] are

$$R_{n0}(x) = \begin{cases} 0, & \text{if } 0 \leq x < \frac{n-1}{2^{k-1}}; \\ \gamma_0 \, 2^{-\frac{(k-1)}{2}-1} \left[C_1(\theta) + C_0(\theta)\right], & \text{if } \frac{n-1}{2^{k-1}} \leq x < \frac{n}{2^{k-1}}; \\ \gamma_0 \, 2^{-\frac{(k-1)}{2}} C_0(\theta), & \text{if } \frac{n}{2^{k-1}} \leq x \leq 1, \end{cases} \quad (2.12)$$

$$R_{n1}(x) = \begin{cases} 0, & \text{if } 0 \leq x < \frac{n-1}{2^{k-1}}; \\ \gamma_1 \, 2^{-\frac{(k-1)}{2}-3} \left[C_2(\theta) - C_0(\theta)\right], & \text{if } \frac{n-1}{2^{k-1}} \leq x < \frac{n}{2^{k-1}}; \\ 0, & \text{if } \frac{n}{2^{k-1}} \leq x \leq 1, \end{cases} \quad (2.13)$$

$$R_{nm}(x) = \begin{cases} 0, & \text{if } 0 \leq x < \frac{n-1}{2^{k-1}}; \\ \gamma_m \, 2^{-\frac{(k-1)}{2}-2} \left[\frac{C_{m+1}(\theta)}{m+1} - \frac{C_{m-1}(\theta)}{m-1} + \mu_m\right], & \text{if } \frac{n-1}{2^{k-1}} \leq x < \frac{n}{2^{k-1}}; \\ \gamma_m \, 2^{-\frac{(k-1)}{2}-2} \rho_m, & \text{if } \frac{n}{2^{k-1}} \leq x \leq 1, \end{cases} \quad (2.14)$$

and

$$S_{n0}(x) = \begin{cases} 0, & \text{if } 0 \leq x < \frac{n-1}{2^{k-1}}; \\ \gamma_0 \, 2^{-\frac{3(k-1)}{2}-4} \left[C_2(\theta) + 4C_1(\theta) + 3C_0(\theta)\right], & \text{if } \frac{n-1}{2^{k-1}} \leq x < \frac{n}{2^{k-1}}; \\ \gamma_0 \, 2^{-\frac{(k-1)}{2}} \left(\frac{1}{2^k} + x - \frac{n}{2^{k-1}}\right), & \text{if } \frac{n}{2^{k-1}} \leq x \leq 1, \end{cases} \quad (2.15)$$

$$S_{n1}(x) = \begin{cases} 0, & \text{if } 0 \leq x < \frac{n-1}{2^{k-1}}; \\ \gamma_1 \, 2^{-\frac{3(k-1)}{2}-4} \left[\frac{C_3(\theta)}{6} - \frac{3C_1(\theta)}{2} - \frac{4C_0(\theta)}{3}\right], & \text{if } \frac{n-1}{2^{k-1}} \leq x < \frac{n}{2^{k-1}}; \\ \gamma_1 \frac{2^{-\frac{3(k-1)}{2}-1}}{-3}, & \text{if } \frac{n}{2^{k-1}} \leq x \leq 1, \end{cases} \quad (2.16)$$

$$S_{n2}(x) = \begin{cases} 0, & \text{if } 0 \leq x < \frac{n-1}{2^{k-1}}; \\ \gamma_2 \, 2^{-\frac{3(k-1)}{2}-3} \left[\frac{C_4(\theta)-1}{24} - \frac{C_2(\theta)-1}{3} - \frac{2}{3}\left(C_1(\theta) + C_0(\theta)\right)\right], & \text{if } \frac{n-1}{2^{k-1}} \leq x < \frac{n}{2^{k-1}}; \\ \gamma_2 \frac{2^{-\frac{(k-1)}{2}}}{-3} \left(\frac{1}{2^k} + x - \frac{n}{2^{k-1}}\right), & \text{if } \frac{n}{2^{k-1}} \leq x \leq 1, \end{cases} \quad (2.17)$$

$$S_{nm}(x) = \begin{cases} 0, & \text{if } 0 \leq x < \frac{n-1}{2^{k-1}}; \\ \gamma_m \, 2^{-\frac{3(k-1)}{2}-3} \begin{bmatrix} \frac{C_{m+2}(\theta)-(-1)^{m+2}}{2(m+1)(m+2)} - \frac{C_m(\theta)-(-1)^m}{2(m+1)(m)} \\ -\frac{C_m(\theta)-(-1)^m}{2(m-1)(m)} + \frac{C_{m-2}(\theta)-(-1)^{m-2}}{2(m-1)(m-2)} \\ + (1 + C_1(\theta))\mu_m \end{bmatrix}, & \text{if } \frac{n-1}{2^{k-1}} \leq x < \frac{n}{2^{k-1}}; \\ \gamma_m \, 2^{-\frac{3(k-1)}{2}-3} \begin{bmatrix} \frac{1-(-1)^{m+2}}{2(m+1)(m+2)} - \frac{1-(-1)^m}{2(m+1)(m)} \\ -\frac{1-(-1)^m}{2(m-1)(m)} + \frac{1-(-1)^{m-2}}{2(m-1)(m-2)} \\ +2\mu_m + 2^k \left(x - \frac{n}{2^{k-1}}\right)\rho_m \end{bmatrix}, & \text{if } \frac{n}{2^{k-1}} \leq x \leq 1, \end{cases} \quad (2.18)$$

where $\theta = 2^k x - 2n + 1$, $\rho_m = \frac{1-(-1)^{m+1}}{m+1} - \frac{1-(-1)^{m-1}}{m-1}$, and $\mu_m = \frac{(-1)^{m-1}}{m-1} - \frac{(-1)^{m+1}}{m+1}$. Moreover, $C_m(x)$ stands for the Chebyshev polynomials of the first kind of degree $m$. For the sake of simplicity and understandability the use of unnecessary notations has been avoided. Thus, for any continuous $f(x)$ throughout this section we have the following relations:

$$f(x) \simeq \sum_{n=1}^{2^{k-1}} \sum_{m=0}^{M-1} d_{nm} I_{nm}(x) = \boldsymbol{D}\boldsymbol{I}(x),$$

$$\int_0^x f(\tau) d\tau \simeq \sum_{n=1}^{2^{k-1}} \sum_{m=0}^{M-1} d_{nm} R_{nm}(x) = \boldsymbol{D}\boldsymbol{R}(x),$$

$$\int_0^x \int_0^\xi f(\tau) d\tau\, d\xi \simeq \sum_{n=1}^{2^{k-1}} \sum_{m=0}^{M-1} d_{nm} S_{nm}(x) = \boldsymbol{D}\boldsymbol{S}(x).$$

## 2.2 Implementation of the method

Let $N_t$ denote the number of divisions of the time interval such that $\Delta t = \frac{T}{N_t}$. On $t \in [t_r, t_{r+1})$ the approximate solution of Eq. (1.1) can be identified as follows:

$$\frac{\partial^3 y}{\partial t \partial x^2}(x,t) \simeq \frac{\partial^3 Y}{\partial t \partial x^2}(x,t) = \sum_{n=1}^{2^{k-1}} \sum_{m=0}^{M-1} d_{nm} I_{nm}(x) = \boldsymbol{D}\boldsymbol{I}(x), \quad t \in [t_r, t_{r+1}). \tag{2.19}$$

Notice that $y(x,t)$ and $Y(x,t)$ represent the exact solution and numerical solution, respectively. Integrating the Eq. (2.19) with respect to $t$ leads to

$$\frac{\partial^2 Y}{\partial x^2}(x,t) = (t - t_r)\boldsymbol{D}\boldsymbol{I}(x) + \frac{\partial^2 Y}{\partial x^2}(x, t_r), \quad t \in [t_r, t_{r+1}). \tag{2.20}$$

On the other hand, integrating Eq. (2.19) twice with respect to $x$ yields

$$\frac{\partial^2 Y}{\partial t \partial x}(x,t) = \boldsymbol{D}\boldsymbol{R}(x) + \frac{\partial^2 Y}{\partial t \partial x}(0,t), \tag{2.21}$$

$$\frac{\partial Y}{\partial t}(x,t) = \boldsymbol{D}\boldsymbol{S}(x) + x\frac{\partial^2 Y}{\partial t \partial x}(0,t) + \frac{\partial Y}{\partial t}(0,t). \tag{2.22}$$

Applying boundary conditions (1.3) in Eq (2.22)

$$\frac{\partial^2 Y}{\partial t \partial x}(0,t) = f_1'(t) - f_0'(t) - \boldsymbol{D}\boldsymbol{S}(1), \tag{2.23}$$

where $'$ stands for the derivative notation. By virtue of the relation given in Eq. (2.23) we have

$$\frac{\partial Y}{\partial t}(x,t) = \boldsymbol{D}\left(\boldsymbol{S}(x) - x \otimes \boldsymbol{S}(1)\right) + x\left(f_1'(t) - f_0'(t)\right) + f_0'(t). \tag{2.24}$$

Additionally, integrating Eq. (2.20) twice with respect to $x$ we get

$$\frac{\partial Y}{\partial x}(x,t) = (t - t_r)\boldsymbol{D}\boldsymbol{R}(x) + \frac{\partial Y}{\partial x}(x, t_r) - \frac{\partial Y}{\partial x}(0, t_r) + \frac{\partial Y}{\partial x}(0, t), \tag{2.25}$$

$$Y(x,t) = (t - t_r)\boldsymbol{D}\boldsymbol{S}(x) + Y(x, t_r) + Y(0, t) - Y(0, t_r)$$
$$+ x\left[\frac{\partial Y}{\partial x}(0, t_r) - \frac{\partial Y}{\partial x}(0, t_r)\right]. \tag{2.26}$$

The required information about $\left[\frac{\partial Y}{\partial x}(0, t_r) - \frac{\partial Y}{\partial x}(0, t_r)\right]$ can be obtained incorporating boundary conditions when $x = 1$. Eq. (2.26) implies

$$\left[\frac{\partial Y}{\partial x}(0, t_r) - \frac{\partial Y}{\partial x}(0, t_r)\right] = f_1(t) - f_0(t) - f_1(t_r) + f_0(t_r) - (t - t_r)\boldsymbol{DS}(1) \quad (2.27)$$

Inserting Eq. (2.27) into both Eq. (2.25) and Eq. (2.26) leads to

$$\frac{\partial Y}{\partial x}(x, t) = (t - t_r)\boldsymbol{D}\left(\boldsymbol{R}(x) - 1 \otimes \boldsymbol{S}(1)\right)$$
$$+ \frac{\partial Y}{\partial x}(x, t_r) + \left(f_1(t) - f_1(t_r) - f_0(t) + f_0(t_r)\right), \quad (2.28)$$

$$Y(x, t) = (t - t_r)\boldsymbol{D}\left(\boldsymbol{S}(x) - x \otimes \boldsymbol{S}(1)\right)$$
$$+ Y(x, t_r) + f_0(t) - f_0(t_r) + x\left(f_1(t) - f_1(t_r) - f_0(t) + f_0(t_r)\right) \quad (2.29)$$

Notice that $\otimes$ denotes the Kronocker product. The required coefficient to generate the approximate solution at $t = t_{r+1}$ can be obtained by substituting the above-mentioned equations into Equation 1.1 such that

$$\boldsymbol{D}\left(\boldsymbol{S}(x) - x \otimes \boldsymbol{S}(1)\right) + x\left(f_1'(t_{r+1}) - f_0'(t_{r+1})\right) + f_0'(t_{r+1}) - A\left[(t_{r+1} - t_r)\boldsymbol{DI}(x) + \frac{\partial^2 Y}{\partial x^2}(x, t_r)\right]$$
$$- B\left[(t - t_r)\boldsymbol{D}\left(\boldsymbol{R}(x) - 1 \otimes \boldsymbol{S}(1)\right) + \frac{\partial Y}{\partial x}(x, t_r) + \left(f_1(t_{r+1}) - f_1(t_r) - f_0(t_{r+1}) + f_0(t_r)\right)\right]$$
$$- X(t_{r+1})\left[(t - t_r)\boldsymbol{D}\left(\boldsymbol{S}(x) - x \otimes \boldsymbol{S}(1)\right) + Y(x, t_r) + f_0(t_{r+1}) - f_0(t_r)\right.$$
$$\left. + x\left(f_1(t_{r+1}) - f_1(t_r) - f_0(t_{r+1}) + f_0(t_r)\right)\right] = \psi(x, t_{r+1}), t \in [t_r, t_{r+1}]. \quad (2.30)$$

It is worth to reminding that by means of the conditions to which the equation subjects we have the following relations at $t = t_0 = 0$.

$$Y(x, 0) = y_0(x) \qquad Y_x(x, 0) = y_0'(x) \qquad Y_{xx}(x, 0) = y_0''(x). \quad (2.31)$$

The source control parameter function needs further investigation in Eq. (2.30) due to the non-availability of the solution. To do so, the overspecified condition is used. Eq. (1.1) can be rewritten at a specific point $x = x_{in}$ as follows:

$$\frac{\partial y}{\partial t}(x_{in}, t) - A\frac{\partial^2 y}{\partial x^2}(x_{in}, t) - B\frac{\partial y}{\partial x}(x_{in}, t) - X(t)y(x_{in}, t) = \psi(x_{in}, t), \ 0 \leq t \leq T \quad (2.32)$$

Using overspecification condition (1.4) in eq. (2.32), we obtain

$$\frac{\partial Q}{\partial t}(t) - A\frac{\partial^2 y}{\partial x^2}(x_{in}, t) - B\frac{\partial y}{\partial x}(x_{in}, t) - X(t)y(x_{in}, t) = \psi(x_{in}, t), \ 0 \leq t \leq T \quad (2.33)$$

In other words,

$$X(t) = \frac{Q'(t) - A\frac{\partial^2 y}{\partial x^2}(x_{in}, t) - B\frac{\partial y}{\partial x}(x_{in}, t) - \psi(x_{in}, t)}{Q(t)}, \ 0 \leq t \leq T. \quad (2.34)$$

By means of the Taylor expansion, at $t = t_{r+1}$, the numerical solution for the source control parameter can be expressed as follows:

$$X(t_{r+1}) = \frac{Q'(t) - A\left[\frac{\partial^2 Y}{\partial x^2}(x_{in}, t_r) + \Delta t \frac{\partial^3 Y}{\partial t \partial x^2}(x_{in}, t_r)\right]}{Q(t)}$$

$$- \frac{B\left[\frac{\partial Y}{\partial x}(x_{in}, t_r) + \Delta t \frac{\partial^2 Y}{\partial t \partial x}(x_{in}, t_r)\right] + \psi(x_{in}, t)}{Q(t)} + \mathcal{O}(\Delta t^2). \quad (2.35)$$

Therefore, Eq. (2.30) can be reduced to an algebraic equation such that

$$M\boldsymbol{D} = \boldsymbol{b}, \quad (2.36)$$

where

$$M = (\boldsymbol{S}(x) - x \otimes \boldsymbol{S}(1)) - A(t_{r+1} - t_r)\boldsymbol{I} - B(t_{r+1} - t_r)(\boldsymbol{R}(x) - 1 \otimes \boldsymbol{S}(1))$$
$$- X(t_{r+1})(t_{r+1} - t_r)(\boldsymbol{S}(x) - x \otimes \boldsymbol{S}(1)), \quad (2.37)$$

and

$$\boldsymbol{b} = \psi(x, t_{r+1}) - x[f_1'(t_{r+1}) - f_0'(t_{r+1})] - f_0'(t_{r+1}) + A\left(\frac{\partial^2 Y}{\partial x^2}(x, t_r)\right)$$
$$+ B\left(\frac{\partial Y}{\partial x}(x, t_r) + f_1(t_{r+1}) - f_1(t_r) - f_0(t_{r+1}) + f_0(t_r)\right)$$
$$+ X(t_{r+1})\left(Y(x, t_r) - f_0(t_r) + f_0(t_{r+1}) + x(f_1(t_{r+1}) - f_1(t_r) - f_0(t_{r+1}) + f_0(t_r))\right). \quad (2.38)$$

After inserting suitable collocation points and using Eq. (2.35), the attained matrix equation Eq. (2.36) is solved at $t = t_{r+1}$ and wavelet coefficients $\boldsymbol{D}$ are obtained for both Taylor and Chebyshev wavelets via the `gmres` package in MATLAB.

# 3 Convergence Analyses

The analysis of the convergence of numerical methods provides a reliability of the method, which enables us to do computational work. For the convergence analysis of the proposed numerical method based on both Taylor wavelet and Chebyshev wavelet, the required functional error estimations of these wavelets should be given. Thereby, this section is designed to state primarily, the error estimations of the mentioned wavelet methods and then to present the convergence result of the proposed method.

## 3.1 The functional accuracy estimations of the Chebyshev wavelet and the Taylor wavelet

### 3.1.1 The error estimations for the Chebyshev wavelet method

**Theorem 3.1.1.** *A function $f(x)$ with bounded second order derivative, that is $|f''(x)| \leq L$, defined on $[0, 1)$, can be expanded as a uniformly convergent infinite series of Chebyshev*

*wavelets as follows:*

$$f(x) = \sum_{n=1}^{\infty} \sum_{m=0}^{\infty} d_{nm} I_{nm}(x).$$

*Notice that $d_{nm} = \langle f(x), I_{nm}(x) \rangle$ where $\langle \cdot, \cdot \rangle$ stands for the inner product in $L_\omega^2[0,1]$ given in Eq. (2.5) and $I_{nm}(x)$ are the Chebyshev wavelets. Moreover,*

$$|d_{nm}| \leq \frac{\gamma_m \pi L}{32 n^{5/2} (m-1)^2}, \quad m > 1. \tag{3.1}$$

*Proof.* The authors refer the interested reader to.[24] □

Let $\sigma_{nm}$ denote the accuracy of the functional approximation of the wavelet method which is defined as

$$\sigma_{nm} = \left\| \int_0^x f(\tau) d\tau - \sum_{n=1}^{2^{k-1}} \sum_{m=0}^{M-1} d_{nm} R_{nm}(x) \right\|.$$

**Theorem 3.1.2.** *A function $f(x)$ with bounded second order derivative, that is $|f''(x)| \leq L$, defined on $[0,1)$. Then, we have the following error estimations for the Chebyshev wavelet method (CWM) and its integral forms:*

$$\sigma_{nm} \leq \frac{\sqrt{\pi} L}{8} \left( \sum_{n=2^{k-1}+1}^{\infty} \sum_{m=M}^{\infty} \frac{1}{n^5} \frac{1}{(m-1)^4} \right)^{1/2} \tag{3.2}$$

*Proof.* The functional accuracy of CWM, i. e. $I_{nm}$, has already given in the study of Sohrabi,.[24] However, for this study a further error analysis is required to investigate since the proposed method includes both $R_{nm}$ and $S_{nm}$ which stands for the first integral and the second integral approaches of the CWM, respectively.

Let us start with the first integral form such that

$$\left\| \int_0^x f(\tau) d\tau - \sum_{n=1}^{2^{k-1}} \sum_{m=0}^{M-1} d_{nm} R_{nm}(x) \right\|^2 = \int_0^1 \left| \sum_{n=1}^{\infty} \sum_{m=0}^{\infty} d_{nm} R_{nm}(x) - \sum_{n=1}^{2^{k-1}} \sum_{m=0}^{M-1} d_{nm} R_{nm}(x) \right|^2 \omega(x) dx$$

This implies that

$$\begin{aligned}
\sigma_{nm}^2 &\leq \int_0^1 \sum_{n=2^{k-1}+1}^{\infty} \sum_{m=M}^{\infty} |d_{nm}|^2 \int_0^x I_{nm}^2(\tau) \omega(\tau) d\tau \omega(x) dx \\
&\leq \sum_{n=2^{k-1}+1}^{\infty} \sum_{m=M}^{\infty} |d_{nm}|^2 \int_0^1 \int_{\frac{n-1}{2^{k-1}}}^{\frac{n}{2^{k-1}}} \gamma_m^2 \, 2^{k-1} \frac{C_m^2 (2^k \tau - 2n + 1)}{\sqrt{1 - (2^k \tau - 2n + 1)^2}} d\tau \omega(x) dx \\
&\leq \sum_{n=2^{k-1}+1}^{\infty} \sum_{m=M}^{\infty} |d_{nm}|^2 \int_0^1 \int_{-1}^{1} \gamma_m^2 \frac{C_m^2(t)}{2\sqrt{1-t^2}} dt \omega(x) dx
\end{aligned}$$

where $t = (2^k\tau - 2n + 1)$. Due to the orthogonality property of the Chebyshev polynomials one can obtain

$$\left\| \int_0^x f(\tau)d\tau - \sum_{n=1}^{2^{k-1}} \sum_{m=0}^{M-1} d_{nm} R_{nm}(x) \right\|^2 \leq \sum_{n=2^{k-1}+1}^{\infty} \sum_{m=M}^{\infty} |d_{nm}|^2 \gamma_m^2 \frac{\pi}{4},$$

By virtue of $|\gamma_m| \leq \frac{2}{\sqrt{\pi}}, \quad m = 0, 1, 2, \ldots$ and Eq. (3.1) we have

$$\sigma_{nm}^2 \leq \frac{\sqrt{\pi}L}{8} \left( \sum_{n=2^{k-1}+1}^{\infty} \sum_{m=M}^{\infty} \frac{1}{n^5(m-1)^4} \right)^{1/2}. \qquad (3.3)$$

The second integral form of the Chebyshev wavelet method follows the line of the first integral form. Besides the use of the bound of $\gamma_m$ and Eq. (3.1) some integration properties has been utilized. This yields

$$\left\| \int_0^x \int_0^u f(\tau)d\tau du - \sum_{n=1}^{2^{k-1}} \sum_{m=0}^{M-1} d_{nm} S_{nm}(x) \right\|^2 \leq \int_0^1 \sum_{n=2^{k-1}+1}^{\infty} \sum_{m=M}^{\infty} |d_{nm}|^2 S_{nm}(x) \omega(x)dx,$$

$$\sigma_{nm}^2 \leq \int_0^1 \sum_{n=2^{k-1}+1}^{\infty} \sum_{m=M}^{\infty} |d_{nm}|^2 \int_0^1 \int_0^u I_{nm}^2(\tau)\omega(\tau)d\tau\omega(u)du\omega(x)dx$$

$$\leq \frac{\sqrt{\pi}L}{8} \left( \sum_{n=2^{k-1}+1}^{\infty} \sum_{m=M}^{\infty} \frac{1}{n^5(m-1)^4} \right)^{1/2}. \qquad (3.4)$$

□

### 3.1.2 The error estimations for the Taylor wavelet method

To the best of our knowledge there is no such analysis for the Taylor wavelet method.

**Theorem 3.1.3.** *A function $f(x)$ with bounded second order derivative, that is $|f''(x)| \leq L$, defined on $[0, 1)$, can be expanded as an infinite sum of Taylor wavelets. Moreover, the series converges uniformly to $f(x)$. That is,*

$$f(x) = \sum_{n=1}^{\infty} \sum_{m=0}^{\infty} d_{nm} I_{nm}(x).$$

*Notice that $d_{nm} = \langle f(x), I_{nm}(x) \rangle$ where $\langle \cdot, \cdot \rangle$ stands for the inner product in $L^2[0,1]$ given in Eq. (2.5) and $I_{nm}(x)$ are the Taylor wavelets. Moreover,*

$$|d_{nm}| \leq \frac{L\sqrt{2m+1}}{n^{\frac{5}{2}}(m+1)(m+2)(m+3)}, \quad m > 1. \qquad (3.5)$$

*Proof.* Using the definition of $d_{nm}$ we have

$$\begin{aligned} d_{nm} &= \int_0^1 f(x) I_{nm}(x) dx, \\ &= \int_{\frac{n-1}{2^{k-1}}}^{\frac{n}{2^{k-1}}} f(x) 2^{\frac{k-1}{2}} \hat{T}_m(2^{k-1}x - n + 1) dx. \end{aligned} \quad (3.6)$$

Changing variables $\tau = 2^k x - n + 1$ leads to

$$d_{nm} = \int_0^1 f\left(\frac{\tau + n - 1}{2^{k-1}}\right) 2^{-\frac{k-1}{2}} \sqrt{2m+1} (\tau)^m d\tau.$$

With the help of the technique of twice integration by parts we get

$$|d_{nm}| \leq \frac{\sqrt{2m+1}}{2^{5\frac{k-1}{2}}(m+1)(m+2)} \left| \int_0^1 f''\left(\frac{\tau + n - 1}{2^{k-1}}\right) (\tau)^{m+2} \right|.$$

One can write by using Cauchy-Schwartz inequality

$$\begin{aligned} |d_{nm}| &\leq \frac{\sqrt{2m+1}}{2^{5\frac{k-1}{2}}(m+1)(m+2)} \left| \int_0^1 f''\left(\frac{\tau + n - 1}{2^{k-1}}\right) \right| \left| \int_0^1 (\tau)^{m+2} \right|, \\ &\leq \frac{L\sqrt{2m+1}}{n^{\frac{5}{2}}(m+1)(m+2)(m+3)}. \end{aligned}$$

As a result of Eq. (3.5), $\sum_{n=1}^{2^{k-1}} \sum_{m=0}^{M-1} d_{nm}$ is absolutely convergent as $k, M \to \infty$. This guarantees the uniform convergence of $\sum_{n=1}^{2^{k-1}} \sum_{m=0}^{M-1} d_{nm} I_{nm}(x)$ to the function $f(x)$. □

In addition to Theorem 3.1.3, Theorem 3.1.4 describes the accuracy of both the Taylor wavelet method and its integral forms. Let $\sigma_{nm}$ stand for the accuracy of the Taylor wavelet method.

**Theorem 3.1.4.** *A function $f(x)$ with bounded second order derivative, that is $|f''(x)| \leq L$, defined on $[0, 1)$. Then, we have the following error estimation for both the Taylor wavelet method (TWM) and its integral forms:*

$$\sigma_{nm} \leq L \left( \sum_{n=2^{k-1}+1}^{\infty} \sum_{m=M}^{\infty} \frac{2m+1}{n^5(m+1)^2(m+2)^2(m+3)^2} \right)^{1/2} \quad (3.7)$$

*Proof.* To estimate the accuracy of the TWM and its integral forms we use the standard technique as follows:

$$\left\| f(x) - \sum_{n=1}^{2^{k-1}} \sum_{m=0}^{M-1} d_{nm} I_{nm}(x) \right\|^2 = \int_0^1 \left| \sum_{n=1}^{\infty} \sum_{m=0}^{\infty} d_{nm} I_{nm}(x) - \sum_{n=1}^{2^{k-1}} \sum_{m=0}^{M-1} d_{nm} I_{nm}(x) \right|^2 dx$$

$$\sigma_{nm}^2 \leq \int_0^1 \sum_{n=2^{k-1}+1}^{\infty} \sum_{m=M}^{\infty} |d_{nm}|^2 I_{nm}^2(x) dx. \quad (3.8)$$

Using the triangle inequality implies that

$$\sigma_{nm}^2 \leq \sum_{n=2^{k-1}+1}^{\infty} \sum_{m=M}^{\infty} |d_{nm}|^2 \int_0^1 I_{nm}^2(x)dx,$$

$$\leq \sum_{n=2^{k-1}+1}^{\infty} \sum_{m=M}^{\infty} |d_{nm}|^2 \int_{\frac{n-1}{2^{k-1}}}^{\frac{n}{2^{k-1}}} \left(2^{\frac{k-1}{2}}\sqrt{2m+1}(2^{k-1}x-n+1)^m\right)^2 dx, \quad (3.9)$$

$$\leq \sum_{n=2^{k-1}+1}^{\infty} \sum_{m=M}^{\infty} |d_{nm}|^2. \quad (3.10)$$

Due to Eq. (3.5) the accuracy estimation of the TWM can be defined as follows:

$$\sigma_{nm} \leq L \left( \sum_{n=2^{k-1}+1}^{\infty} \sum_{m=M}^{\infty} \frac{2m+1}{n^5(m+1)^2(m+2)^2(m+3)^2} \right)^{1/2}$$

On the other hand, with the help of the orthogonality property the accuracy of the integral forms of the TWM can be obtain as follows:

$$\left\| \int_0^x f(\tau)d\tau - \sum_{n=1}^{2^{k-1}} \sum_{m=0}^{M-1} d_{nm}R_{nm}(x) \right\|^2 = \int_0^1 \left| \sum_{n=1}^{\infty} \sum_{m=0}^{\infty} d_{nm}R_{nm}(x) - \sum_{n=1}^{2^{k-1}} \sum_{m=0}^{M-1} d_{nm}R_{nm}(x) \right|^2 dx.$$

This implies that

$$\sigma_{nm}^2 \leq \sum_{n=2^{k-1}+1}^{\infty} \sum_{m=M}^{\infty} |d_{nm}|^2 \int_0^1 \int_0^x I_{nm}^2(\tau) d\tau dx$$

$$\sigma_{nm} \leq L \left( \sum_{n=2^{k-1}+1}^{\infty} \sum_{m=M}^{\infty} \frac{2m+1}{n^5(m+1)^2(m+2)^2(m+3)^2} \right)^{1/2}. \quad (3.11)$$

where $t = \left(2^{k-1}\tau - n + 1\right)$. Likewise, the second integral form of the TWM is bounded by

$$\left\| \int_0^x \int_0^\tau f(u)dud\tau - \sum_{n=1}^{2^{k-1}} \sum_{m=0}^{M-1} d_{nm}S_{nm}(x) \right\|^2 = \int_0^1 \left| \sum_{n=1}^{\infty} \sum_{m=0}^{\infty} d_{nm}S_{nm}(x) - \sum_{n=1}^{2^{k-1}} \sum_{m=0}^{M-1} d_{nm}S_{nm}(x) \right|^2 dx$$

$$\sigma_{nm}^2 \leq \int_0^1 \sum_{n=2^{k-1}+1}^{\infty} \sum_{m=M}^{\infty} |d_{nm}|^2 \int_0^x \int_0^\tau I_{nm}^2(u)dud\tau dx$$

$$\sigma_{nm} \leq L \left( \sum_{n=2^{k-1}+1}^{\infty} \sum_{m=M}^{\infty} \frac{2m+1}{n^5(m+1)^2(m+2)^2(m+3)^2} \right)^{1/2}. \quad (3.12)$$

□

## 3.2 The final convergence result of the proposed method

Under the lights of above-mentioned error estimations section 3.1 is ending up with Theorem 3.2.1 which declares the convergence result of the proposed numerical scheme.

**Theorem 3.2.1.** *Let $y \in L^2[0,1] \cap \mathcal{C}[0,T]$ be the exact solution of initial-boundary problem stated in Eq. (1.1)-(1.4). Let $\Delta t = \frac{T}{N_t}$ where $N_t$ stands for the number of discretizations of the time interval. Suppose that $Y(x,t)$ is the numerical solution obtained by the proposed method. The proposed method is convergent in the sense that*

$$\|y(x,t_r) - Y(x,t_r)\| \leq \|y(x,0) - Y(x,0)\| + r\kappa\Delta t.$$

*Notice that $\kappa = \lambda \sigma_{nm}$, $\lambda \in \mathbb{R}$ where $\sigma_{nm}$ represents the error estimations of the considered wavelet methods. Moreover, $r$ represents the time step for $r = 1, 2, ..., N_t$*

*Proof.* Prior to undertaking the analysis, recall the numerical solution and the exact solution at $t = t_{r+1}$, respectively:

$$Y(x, t_{r+1}) = (t_{r+1} - t_r) \sum_{n=1}^{2^{k-1}} \sum_{m=0}^{M-1} d_{nm} (S_{nm}(x) - x \otimes S_{nm}(1)) + Y(x, t_r) + bc(t_r^{r+1}, x), \tag{3.13}$$

and

$$y(x, t_{r+1}) = (t_{r+1} - t_r) \sum_{n=1}^{\infty} \sum_{m=0}^{\infty} d_{nm} (S_{nm}(x) - x \otimes S_{nm}(1)) + Y(x, t_r) + bc(t_r^{r+1}, x), \tag{3.14}$$

where $bc(t_r^{r+1}, x) = f_0(t_{r+1}) - f_0(t_r) + x(f_1(t_{r+1}) - f_1(t_r) - f_0(t_{r+1}) + f_0(t_r))$. □

Subtracting Eq. (3.13) from Eq. (3.14) the local error can be defined as follows:

$$|Y(x, t_{r+1}) - y(x, t_{r+1})| \leq \Delta t \sum_{n=2^{k-1}+1}^{\infty} \sum_{M}^{\infty} |d_{nm}(S_{nm}(x) - x \otimes S_{nm}(1))| + |Y(x, t_r) - y(x, t_r)|, \tag{3.15}$$

where $\Delta t = (t_{r+1} - t_r)$. Notice that the terms of $bc(t_r^{r+1}, x)$ are extracted from the exact solution, therefore, there is no contribution of $bc(t_r^{r+1}, x)$ on the error estimation. Moreover, Eq. (3.15) is hinted at a deep connection between the convergence result of the proposed method and error estimations of wavelet methods. By defining $e_j = |Y(x, t_j) - y(x, t_{j-1})|$, $j = 1, 2, ..., r+1$, we have

$$e_r \leq e_{r-1} + \kappa \Delta t, \tag{3.16}$$

where

$$\kappa = \begin{cases} \lambda \left( \sum_{n=2^{k-1}+1}^{\infty} \sum_{m=M}^{\infty} \frac{2m+1}{n^5(m+1)^2(m+2)^2(m+3)^2} \right)^{1/2}, & \text{for TWM}; \\ \lambda \left( \sum_{n=2^{k-1}+1}^{\infty} \sum_{m=M}^{\infty} \frac{1}{n^5(m-1)^4} \right)^{1/2}, & \text{for CWM}. \end{cases} \tag{3.17}$$

Here, $\lambda$ is a constant depends on $L$ and $\max_{1 \leq i \leq 2^{k-1}M} x_i$. The error propagation can easily seen by induction as follows:

$$e_1 \leq e_0 + \kappa \Delta t \tag{3.18}$$
$$e_2 \leq e_1 + \kappa \Delta t \leq e_0 + 2\kappa \Delta t \tag{3.19}$$
$$\vdots \quad \vdots \quad \vdots \tag{3.20}$$
$$e_{N_t} \leq e_0 + N_t \kappa \Delta t, \tag{3.21}$$

Moreover, it is worth noting that at $t = 0$ the numerical solution is obtained by the initial condition of the equation, that is $e_0 = 0$. Furthermore, based on the definition of $\Delta t$ one can say that $N_t \Delta t = T$. Therefore, by virtue of $\sigma_{nm} \to 0$ as $k$ and $M$ increases one can be concluded that the proposed method is convergent by $e_{N_t} \to 0$.

## 4 Numerical Results and Discussion

In this section, we present numerical results based on Taylor wavelet method and Chebyshev wavelet method applied on one dimensional parabolic inverse problem. Numerical results are in good agreement with the exact results. To show the accuracy of the method, we have provided $L_\infty$ and $L_2$ error of the solution. The $L_\infty$ and $L_2$ error norm are defined as follows:

$$\|y(.,t_r) - Y(.,t_r)\|_{L_\infty} = \max_{1 \leq l \leq 2^{k-1}M} |y(x_l, t_r) - Y(x_l, t_r)|, \tag{4.1}$$

$$\|y(.,t_r) - Y(.,t_r)\|_{L_2} = \frac{1}{2^{k-1}M} \left( \sum_{l=1}^{2^{k-1}M} |y(x_l, t_r) - Y(x_l, t_r)|^2 \right)^{1/2}. \tag{4.2}$$

In the numerical results CWM denotes Chebyshev wavelet method whereas TWM denotes Taylor wavelet method.

**Example 1.**

$$\frac{\partial y}{\partial t}(x,t) = \frac{\partial^2 y}{\partial x^2}(x,t) + 2\frac{\partial y}{\partial x}(x,t) + X(t)y(x,t) - (2 + xt^2)e^t, \ 0 \leq x \leq 1, 0 < t \leq T, \tag{4.3}$$

with initial condition

$$y(x,0) = x, \qquad 0 \leq x \leq 1, \tag{4.4}$$

and Dirichlet boundary conditions

$$y(0,t) = 0, \qquad 0 < t \leq T, \tag{4.5}$$

$$y(1,t) = e^t, \qquad 0 < t \leq T, \tag{4.6}$$

subject to the overspecified condition

$$y(0.5, t) = \frac{e^t}{2}, \qquad 0 < t \leq T. \tag{4.7}$$

The exact solution is

$$y(x,t) = xe^t,$$

and
$$X(t) = 1 + t^2.$$

| $x$ | Pointwise absolute error | | | |
| --- | --- | --- | --- | --- |
| | CWM | TWM | CWM | TWM |
| | $\Delta t = 10^{-2}$ | $\Delta t = 10^{-2}$ | $\Delta t = 10^{-3}$ | $\Delta t = 10^{-3}$ |
| 0.125 | 9.812 $\times 10^{-3}$ | 5.982 $\times 10^{-5}$ | 3.702 $\times 10^{-4}$ | 4.062 $\times 10^{-7}$ |
| 0.250 | 1.738 $\times 10^{-2}$ | 1.012 $\times 10^{-4}$ | 6.106 $\times 10^{-4}$ | 6.784 $\times 10^{-7}$ |
| 0.375 | 2.245 $\times 10^{-2}$ | 1.252 $\times 10^{-4}$ | 8.021 $\times 10^{-4}$ | 8.303 $\times 10^{-7}$ |
| 0.500 | 2.452 $\times 10^{-2}$ | 1.329 $\times 10^{-4}$ | 8.834 $\times 10^{-4}$ | 8.721 $\times 10^{-7}$ |
| 0.625 | 2.318 $\times 10^{-2}$ | 1.220 $\times 10^{-4}$ | 9.265 $\times 10^{-4}$ | 8.115 $\times 10^{-7}$ |
| 0.750 | 1.854 $\times 10^{-2}$ | 9.709 $\times 10^{-5}$ | 7.580 $\times 10^{-4}$ | 6.657 $\times 10^{-7}$ |
| 0.875 | 1.073 $\times 10^{-2}$ | 5.632 $\times 10^{-5}$ | 4.562 $\times 10^{-4}$ | 4.356 $\times 10^{-7}$ |

Table 1: Comparison results for the solution at $t = 1.0, k = 4$ and $M = 4$.

| $t$ | Exact $X$ | Pointwise absolute error | |
| --- | --- | --- | --- |
| | | CWM | TWM |
| 0.1 | 1.01 | $3.545 \times 10^{-5}$ | $3.773 \times 10^{-6}$ |
| 0.2 | 1.04 | $1.547 \times 10^{-5}$ | $3.434 \times 10^{-6}$ |
| 0.3 | 1.09 | $3.717 \times 10^{-4}$ | $3.227 \times 10^{-6}$ |
| 0.4 | 1.16 | $1.125 \times 10^{-3}$ | $3.122 \times 10^{-6}$ |
| 0.5 | 1.25 | $2.632 \times 10^{-3}$ | $3.085 \times 10^{-6}$ |
| 0.6 | 1.36 | $4.494 \times 10^{-3}$ | $3.096 \times 10^{-6}$ |
| 0.7 | 1.49 | $6.549 \times 10^{-3}$ | $3.144 \times 10^{-6}$ |
| 0.8 | 1.64 | $4.703 \times 10^{-3}$ | $3.229 \times 10^{-6}$ |
| 0.9 | 1.81 | $3.523 \times 10^{-3}$ | $3.361 \times 10^{-6}$ |
| 1.0 | 2.00 | $3.043 \times 10^{-3}$ | $3.557 \times 10^{-6}$ |

Table 2: Comparison results for the control parameter at $\Delta t = 10^{-3}, k = 4$ and $M = 4$.

In Table 1, we have presented pointwise absolute error of solutions obtained by Chebyshev wavelet method and Taylor wavelet method at different $\Delta t$ whereas Table 2 presents pointwise absolute error of the control parameter at $\Delta t = 10^{-3}$. It is evident from the numerical results that Taylor wavelet method performs better than the Chebyshev wavelet method. As expected, when we are reducing the time step size, pointwise absolute errors are decaying fastly. It is also observed that at very less $k$ and $M$, one can obtain a very good accuracy. The CPU time taken by Taylor wavelet method is 0.15s and 0.60s when $\Delta t = 10^{-2}$ and $10^{-3}$ whereas CPU time taken by Chebyshev wavelet method is 0.18s and 0.80s respectively when $\Delta t = 10^{-2}$ and $10^{-3}$.

Figure 1 and Figure 2 present $L_\infty$ and $L_2$ error of the solution while solving the parabolic inverse problem using Chebyshev wavelet method and Taylor wavelet method respectively. It is observed that both the errors decay fastly with increasing time.

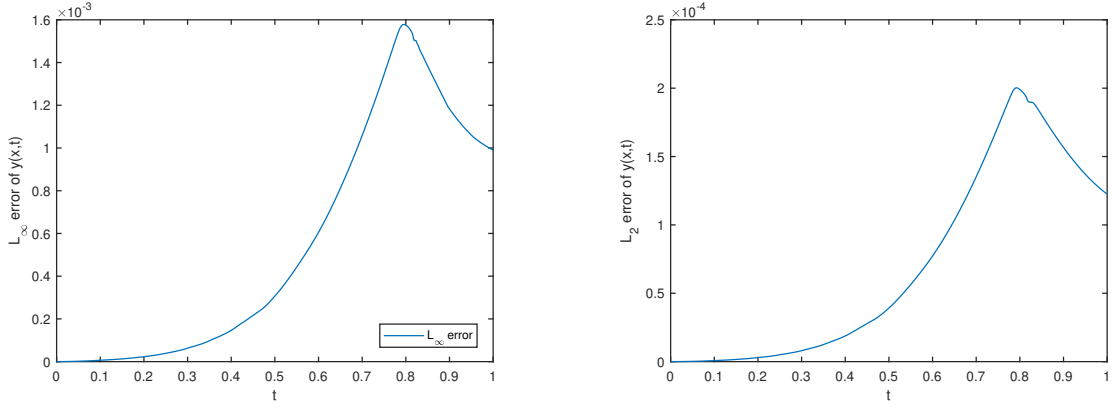

Figure 1: $L_\infty$ and $L_2$ errors of the Chebyshev wavelet solution at $\Delta t = 10^{-3}, k = 4$ and $M = 4$.

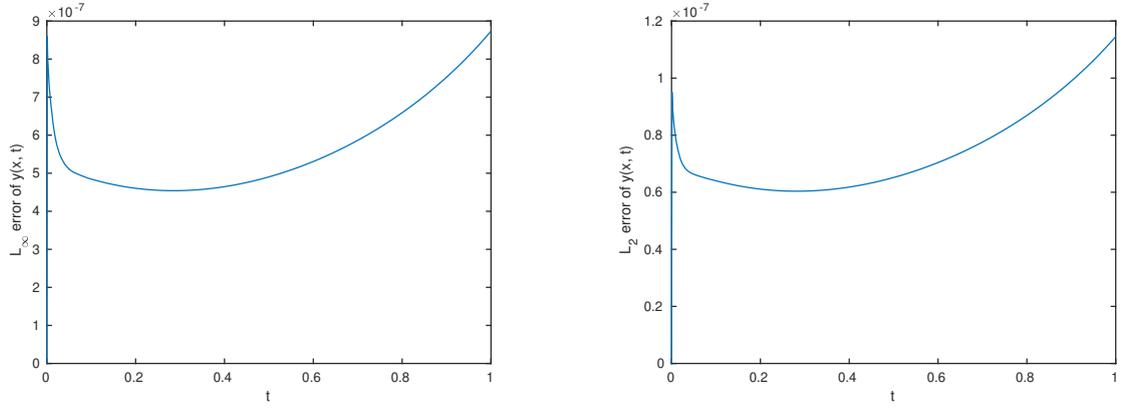

Figure 2: $L_\infty$ and $L_2$ errors of Taylor wavelet solution at $\Delta t = 10^{-3}, k = 4$ and $M = 4$.

**Example 2.**

$$\frac{\partial y}{\partial t} = \frac{\partial^2 y}{\partial x^2} + X(t)y + x\cos(t) - tx\sin(t), \qquad 0 \leq x \leq 1, 0 < t \leq T, \qquad (4.8)$$

with initial condition
$$y(x, 0) = 0, \qquad 0 \leq x \leq 1, \qquad (4.9)$$

and Dirichlet boundary conditions
$$y(0, t) = 0, \qquad 0 < t \leq T, \qquad (4.10)$$
$$y(1, t) = \sin(t), \qquad 0 < t \leq T, \qquad (4.11)$$

subject to the overspecified condition
$$y(0.5, t) = 0.5 \sin t, \qquad 0 < t \leq T. \qquad (4.12)$$

The exact solution is
$$y(x, t) = x \sin t,$$

and
$$X(t) = t.$$

| $x$ | Pointwise absolute error | | | |
| | CWM | TWM | CWM | TWM |
| | $\Delta t = 10^{-2}$ | $\Delta t = 10^{-2}$ | $\Delta t = 10^{-3}$ | $\Delta t = 10^{-3}$ |
| 0.125 | $3.449 \times 10^{-4}$ | $1.500 \times 10^{-5}$ | $1.390 \times 10^{-5}$ | $3.369 \times 10^{-8}$ |
| 0.250 | $6.737 \times 10^{-4}$ | $2.850 \times 10^{-5}$ | $2.845 \times 10^{-5}$ | $5.933 \times 10^{-8}$ |
| 0.375 | $9.652 \times 10^{-4}$ | $3.895 \times 10^{-5}$ | $4.241 \times 10^{-5}$ | $8.526 \times 10^{-8}$ |
| 0.500 | $1.168 \times 10^{-3}$ | $4.503 \times 10^{-5}$ | $5.341 \times 10^{-5}$ | $9.583 \times 10^{-8}$ |
| 0.625 | $1.230 \times 10^{-3}$ | $4.529 \times 10^{-5}$ | $5.772 \times 10^{-5}$ | $1.043 \times 10^{-7}$ |
| 0.750 | $1.093 \times 10^{-3}$ | $3.851 \times 10^{-5}$ | $5.152 \times 10^{-5}$ | $1.022 \times 10^{-7}$ |
| 0.875 | $6.955 \times 10^{-4}$ | $2.342 \times 10^{-5}$ | $3.183 \times 10^{-5}$ | $7.392 \times 10^{-8}$ |

Table 3: Comparison results for the solution at $t = 0.5, k = 4$ and $M = 4$.

| $t$ | Exact $X$ | Pointwise absolute error | |
| | | CWM | TWM |
| 0.05 | 0.05 | $1.248 \times 10^{-4}$ | $5.833 \times 10^{-5}$ |
| 0.10 | 0.10 | $2.022 \times 10^{-4}$ | $2.385 \times 10^{-5}$ |
| 0.15 | 0.15 | $1.516 \times 10^{-4}$ | $1.283 \times 10^{-5}$ |
| 0.20 | 0.20 | $4.838 \times 10^{-5}$ | $7.842 \times 10^{-6}$ |
| 0.25 | 0.25 | $6.991 \times 10^{-5}$ | $5.373 \times 10^{-6}$ |
| 0.30 | 0.30 | $1.968 \times 10^{-4}$ | $3.869 \times 10^{-6}$ |
| 0.35 | 0.35 | $2.988 \times 10^{-4}$ | $2.888 \times 10^{-6}$ |
| 0.40 | 0.40 | $3.246 \times 10^{-4}$ | $2.219 \times 10^{-6}$ |
| 0.45 | 0.45 | $3.199 \times 10^{-4}$ | $1.749 \times 10^{-6}$ |
| 0.50 | 0.50 | $3.397 \times 10^{-4}$ | $1.409 \times 10^{-6}$ |

Table 4: Comparison results for the control parameter at $\Delta t = 10^{-3}, k = 4$ and $M = 4$.

In Table 3, we have presented pointwise absolute error of solutions obtained by Chebyshev wavelet method and Taylor wavelet method at different $\Delta t$ whereas Table 4 presents pointwise absolute error of the control parameter at $\Delta t = 10^{-3}$. The CPU time taken by Taylor wavelet method is 0.14s and 0.64s when $\Delta t = 10^{-2}$ and $10^{-3}$ whereas CPU time taken by Chebyshev wavelet method is 0.15s and 0.70s respectively when $\Delta t = 10^{-2}$ and $10^{-3}$.

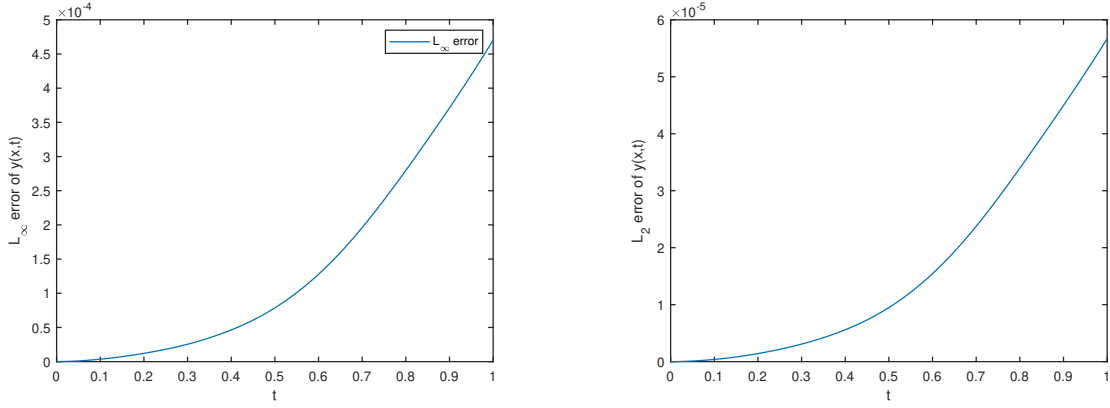

Figure 3: $L_\infty$ and $L_2$ errors of Chebyshev wavelet solution at $\Delta t = 10^{-3}, k = 4$ and $M = 4$.

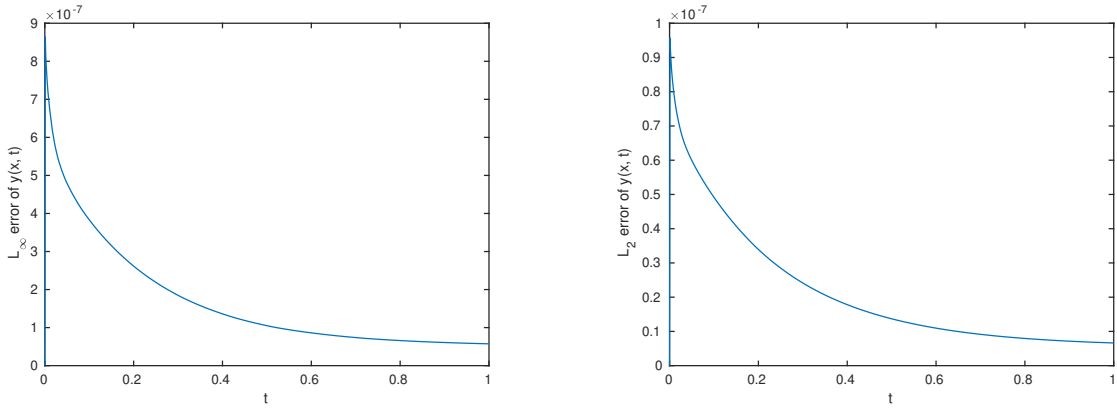

Figure 4: $L_\infty$ and $L_2$ errors of Taylor wavelet solution at $\Delta t = 10^{-3}, k = 4$ and $M = 4$.

Figure 3 and Figure 4 present $L_\infty$ and $L_2$ error of the solution while solving the parabolic inverse problem using Chebyshev wavelet method and Taylor wavelet method respectively.

# 5 Conclusion

We have developed efficient and accurate numerical methods based on Taylor and Chebyshev wavelets for the parameter identification in parabolic inverse problem. Using the uniform convergence property possessed by Taylor and Chebyshev wavelet, we have derived a rigorous convergence analysis. We have compared both the proposed methods on two parabolic inverse problem and it is evident from the numerical result that Taylor wavelet based method provides better result than the Chebyshev wavelet method. Considering very few collocation points in the domain, we are able to achieve a very good accuracy. Due to very less computational cost, CPU time taken to solve the problem is very less. The proposed method can be extended to higher dimensional problems easily.